\documentclass{elsarticle}  %
\usepackage{graphicx}
\usepackage{amsmath}
\usepackage{amssymb}
\usepackage{enumerate}
\usepackage{theorem}
\usepackage{rawfonts}
\usepackage{setspace}

\usepackage{arydshln}
\usepackage{lscape}

\usepackage{float}
\usepackage{color}
\usepackage{multicol}

\usepackage[all]{xy}
\xyoption{poly} \xyoption{graph} \xyoption{arc} \xyoption{web}
\xyoption{2cell}

\begingroup
\newtheorem{thrm}{Theorem}[section]

\newtheorem{lemm}[thrm]{Lemma\hskip 1mm}
\newtheorem{exa} [thrm]{Example\hskip 1mm}
\newtheorem{dfn} [thrm]{Definition}

\endgroup

\newenvironment{pf}{\noindent\textbf{Proof}}{\hspace*{\fill}$\square$\\[6pt]}

\begin{document}

\title{On connected degree sequences}

\author{Jonathan McLaughlin}
\address{Department of Mathematics, St. Patrick's College, Dublin City University, Dublin 9, Ireland }
\ead{jonny$\_\:$mclaughlin@hotmail.com}

\parindent=0cm


\begin{abstract} This note gives necessary and sufficient conditions for a sequence of non-negative integers to be the degree sequence of a connected simple graph. This result is implicit in a paper of Hakimi. A new alternative characterisation of these necessary and sufficient conditions is also given. 

\end{abstract}

\begin{keyword} connected graph \sep degree sequence 
\MSC[2010] 05C40
\end{keyword}

\maketitle

\parindent=0cm

\section{Introduction}
A finite sequence of non-negative integers is called a graphic sequence if it is the degree sequence of some finite simple graph. Erd\"os and Gallai \cite{EG} first found necessary and sufficient conditions for a sequence of non-negative integers to be graphic and these conditions have since been refined by Hakimi \cite{Hk} (stated in the sequel as Theorem \ref{HH}) as well as (independently) by Havel \cite{Hv}. Alternative characterisations and generalisations are due to Choudum \cite{Ch}, Sierksma $\And$ Hoogeveen \cite{SH} and Tripathi et al. \cite{TVW10}, \cite{TV03}, \cite{TV07}. This note states a result which is implicit in Hakimi \cite{Hk}, before giving an alternative characterisation of these necessary and sufficient conditions for a finite sequence of non-negative integers to be the degree sequence of a connected simple graph.

\section{Preliminaries }

Let $G=(V_{G},E_{G})$ be a graph where $V_{G}$ denotes the vertex set of $G$ and $E_{G}\subseteq [V_{G}]^{2}$ denotes the edge set of $G$ (given that $[V_G]^2$ is the set of all $2$-element subsets of $V_G$).  An edge $\{a,b\}$ is denoted $ab$ in the sequel. A graph is {\it finite} when $|V_{G}|<\infty$ and $|E_{G}|<\infty$, where $|X|$ denotes the cardinality of the set $X$. A graph is {\it simple} if it contain no loops (i.e. $a\neq b$ for all $ab\in E_{G})$ or parallel/multiple edges (i.e. $E_{G}$ is not a multiset). The {\it degree} of a vertex $v$ in a graph $G$, denoted $deg(v)$, is the number of edges in $G$ which contain $v$. A {\it path} is a graph with $n$ vertices in which two vertices, known as the {\it endpoints}, have degree $1$ and $n-2$ vertices have degree $2$. A graph is {\it connected} if there exists at least one path between every pair of vertices in the graph. A {\it tree} is a connected graph with $n$ vertices and $n-1$ edges. $K_{n}$ denotes the {\it complete graph} on $n$ vertices. All basic graph theoretic definitions can be found in standard texts such as \cite{BM}, \cite{D} or \cite{GG}. All graphs in this note are undirected and finite. 


\section{Degree sequences and graphs }\label{s4}

A finite sequence $s=\{s_1,...,s_n\}$ of non-negative integers is called {\it realisable} if there exists a finite graph with vertex set  $\{v_1,..., v_n\}$ such that $deg(v_i)=s_i$ for all $i=1,...,n$. A sequence $s$ which is realisable as a simple graph is called {\it graphic}. Given a graph $G$ then the {\it degree sequence} of $G$, denoted $d(G)$, is the monotonic non-increasing sequence of degrees of the vertices in $V_G$. This means that every realisable (resp. graphic) sequence $s$ is equal to the degree sequence $d(G)$ of some graph (resp. simple graph) $G$ (subject to possible rearrangement of the terms in $s$). The maximum degree of a vertex in $G$ is denoted $\Delta_G$ and the minimum degree of a vertex in $G$ is denoted $\delta_G$. In this note all sequences will have {\it positive} terms as the only connected graph which has a degree sequence containing a zero is $(\{v\},\{\})$.  \\

The following theorem states necessary and sufficient conditions for a sequence to be realisable (though not necessarily graphic).

\begin{thrm}[Hakimi]\label{H2} Given a sequence $s=\{s_1,...,s_n\}$ of positive integers such that $s_i\geq s_{i+1}$ for $i=1,...,n-1$ then $s$ is realisable if and only if $\sum\limits_{i=1}^{n}s_i $ is even and $\sum\limits_{i=2}^{n}s_i \geq s_1 $.
\end{thrm}

To address the issue of when a sequence is graphic, Hakimi describes in \cite{Hk} a process he called a {\it reduction cycle} and uses it to state the following result.  
\begin{thrm}[Havel, Hakimi]\label{HH} Given a sequence $\{s_1,...,s_n\}$ of positive integers such that $s_i\geq s_{i+1}$ for $i=1,...,n-1$ then the sequence $\{s_1,...,s_n\}$ is graphic if and only if the sequence $\{s_2-1,s_3-1,...,s_{s_1 +1}-1, s_{s_1+2},...,s_n\}$ is graphic.
\end{thrm}

\section{Degree sequences and connected graphs }\label{s5}

\begin{dfn} A finite sequence $s=\{s_1,...,s_n\}$ of positive integers is called {\it connected} (resp. connected and graphic) if $s$ is realisable as a connected graph (resp. connected simple graph) with vertex set  $\{v_1,..., v_n\}$ such that $v_i$ has degree $s_i$ for all $i=1,...,n$.
\end{dfn}

Of course disconnected realisations of connected and graphic degree sequences exist, for example, $(2,2,2,2,2,2)$ can be realised as a $6$-cycle or as two disjoint $3$-cycles.\\ 

As a graph is connected if and only if it contains a spanning tree, then a simple induction argument on the number of edges shows that every spanning tree of a graph $G$, with $|V_G|=n$, has exactly $n-1$ edges. Hence, a necessary condition for a graph $G$, with $|V_G|=n$, to be connected is that $|E_G|\geq n-1$.  \\

The following theorem states necessary and sufficient conditions for a sequence to be connected but not necessarily simple.

\begin{thrm}[Hakimi]\label{H1} Given a sequence $s=\{s_1,...,s_n\}$ of positive integers such that $s_i\geq s_{i+1}$ for $i=1,...,n-1$ then $s$ is connected if and only if $s$ is realisable and $\sum\limits_{i=1}^{n}s_i \geq 2(n-1)$.
\end{thrm}

The two main tools used in the proof of Theorem \ref{H1} are $d${\it-invariant operations} (which leave degree sequences unchanged) and Lemma \ref{LH} (which appears in \cite{Hk} as Lemma $1$). \\

\begin{figure}[H]
\begin{center}
\scalebox{0.75}{$\begin{xy}\POS (0,5) *\cir<2pt>{} ="a" *+!DR{a},
(15,5) *\cir<2pt>{} ="b" *+!DL{b},
 (0,-5) *\cir<2pt>{} ="c" *+!UR{c},
(15,-5) *\cir<2pt>{} ="d" *+!UL{d},
  (-10,12)*+!{G},
  
\POS "a" \ar@{-}  "b",
\POS "c" \ar@{-}  "d",

\POS(7.5,0),  {\ellipse(20,12)<>{}},

\POS  (35,0)*+!{\longrightarrow},

\POS (55,5) *\cir<2pt>{} ="a" *+!DR{a},
(70,5) *\cir<2pt>{} ="b" *+!DL{b},
 (55,-5) *\cir<2pt>{} ="c" *+!UR{c},
(70,-5) *\cir<2pt>{} ="d" *+!UL{d},
  (45,12)*+!{G'},
  
\POS "a" \ar@{-}  "c",
\POS "b" \ar@{-}  "d",

\POS(62.5,0),  {\ellipse(20,12)<>{}},

\POS  (90,0)*+!{or},

\POS (110,5) *\cir<2pt>{} ="a" *+!DR{a},
(125,5) *\cir<2pt>{} ="b" *+!DL{b},
 (110,-5) *\cir<2pt>{} ="c" *+!UR{c},
(125,-5) *\cir<2pt>{} ="d" *+!UL{d},
  (100,12)*+!{G'},
  
\POS "a" \ar@{-}  "d",
\POS "c" \ar@{-}  "b",

\POS(117.5,0),  {\ellipse(20,12)<>{}},

 \end{xy}$}

\caption{The two possible $d$-invariant operations on $G$ resulting in $G'$ such that $d(G)=d(G')$ }
\label{Cases1}
\end{center}
\end{figure}
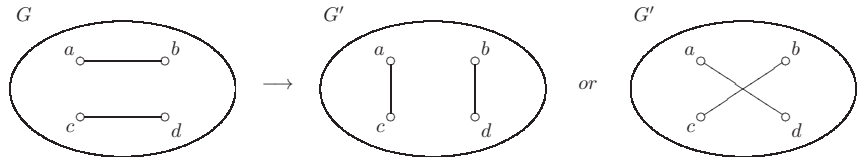

Consider a graph $G$ with $d(G)=(d_1,...,d_n)$. Given any two edges $ab, cd \in E_G$, where $a,b,c$ and $d$ are all distinct, then $G$ is transformed by a $d${\it -invariant operation} into $G'$ when either
\begin{itemize}
\item $V_{G'}=V_{G}$ and $E_{G'}=(E_{G}\setminus\{ab,cd\})\cup \{ac,bd\}$, or  
\item  $V_{G'}=V_{G}$ and $E_{G'}=(E_{G}\setminus\{ab,cd\})\cup \{ad,bc\}$. 
\end{itemize}
Figure \ref{Cases1} shows both $d$-invariant operations. 

These $d$-invariant operations are used to prove the following important result.   

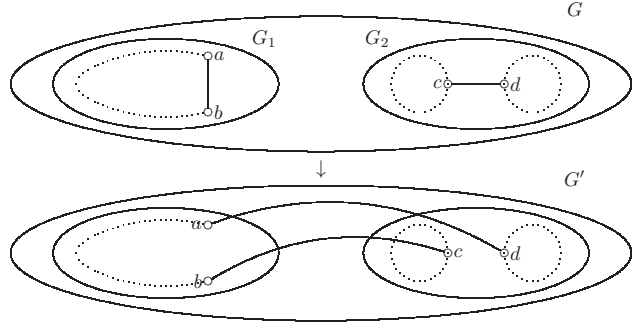
\begin{figure}[H]
\begin{center}
\scalebox{0.75}{$\begin{xy}\POS (0,5) *\cir<0pt>{} ="c" *+!UR{},
(15,5) *\cir<2pt>{} ="a" *+!L{a},
 (0,-5) *\cir<0pt>{} ="d" *+!UR{},
(15,-5) *\cir<2pt>{} ="b" *+!L{b},
  (25,7)*+!{G_1},
  
\POS "a" \ar@{-}  "b",
\POS "a" \ar@{.}@/_0.2pc/  "c",
\POS "d" \ar@{.}@/_0.2pc/  "b",
\POS "c" \ar@{.}@/_2pc/  "d",

\POS(7.5,0),  {\ellipse(20,8)<>{}},

\POS (57.5,0) *\cir<2pt>{} ="c" *+!R{c},
(67.5,0) *\cir<2pt>{} ="d" *+!L{d},
  (80,12)*+!{G},
  (45,7)*+!{G_2},
\POS "c" \ar@{-}  "d",

\POS(62.5,0),  {\ellipse(20,8)<>{}},
\POS(52.5,0),  {\ellipse(5,5)<>{.}},
\POS(72.5,0),  {\ellipse(5,5)<>{.}},

\POS(35,0),  {\ellipse(55,12)<>{}},

\POS  (35,-16)*+!{\downarrow},

\POS(35,-30),  {\ellipse(55,12)<>{}},

\POS  (0,-25) *\cir<0pt>{} ="c" *+!UR{},
(15,-25) *\cir<2pt>{} ="a" *+!R{a},
 (0,-35) *\cir<0pt>{} ="d" *+!UR{},
(15,-35) *\cir<2pt>{} ="b" *+!R{b},
  (57.5,-30) *\cir<2pt>{} ="e" *+!L{c},
(67.5,-30) *\cir<2pt>{} ="f" *+!L{d},
   (80,-18)*+!{G'},
  
\POS "a" \ar@{.}@/_0.2pc/  "c",
\POS "d" \ar@{.}@/_0.2pc/  "b",
\POS "c" \ar@{.}@/_2pc/  "d",

\POS "a" \ar@{-}@/^1.5pc/  "f",
\POS "b" \ar@{-}@/^1.2pc/  "e",

\POS(7.5,-30),  {\ellipse(20,8)<>{}},
\POS(62.5,-30),  {\ellipse(20,8)<>{}},
\POS(52.5,-30),  {\ellipse(5,5)<>{.}},
\POS(72.5,-30),  {\ellipse(5,5)<>{.}},

 \end{xy}$}

\caption{$G'$ has one less connected subgraph than $G$ following a $d$-invariant operation}
\label{Cases2}
\end{center}
\end{figure}

\begin{lemm}\label{LH} Let $G_1,G_2,...,G_r$, (with $r>1$), be maximally connected subgraphs of $G$ such that not all of the $G_i$ are acyclic, then there exists a graph $G'$ with $r-1$ maximally connected subgraphs such that $d(G')=d(G)$.   
\end{lemm}
The essence of Lemma \ref{LH} is presented in Figure \ref{Cases2}. Note that {\it worst-case-scenarios} are assumed i.e. $G_1$ is a cycle and $G_2$ is acyclic with the graph $\big(V_{G_2}, (E_{G_2}\setminus \{cd\})\big)$ being disconnected.

\section{Results}
The first result, Theorem \ref{Implicit}, is an explicit statement of a result implicit in \cite{Hk}. The second result, Theorem \ref{Main}, is a new alternative characterisation of Theorem \ref{Implicit} and has a similar flavour to that of Theorem \ref{HH}. 

\begin{thrm}\label{Implicit} Given a sequence $s=\{s_1,...,s_n\}$ of positive integers such that $s_i\geq s_{i+1}$ for $i=1,...,n-1$ then $s$ is connected and graphic if and only if the sequence $s'=\{s_1',...,s_{n-1}'\}=\{s_2-1,s_3-1,...,s_{s_1 +1}-1, s_{s_1+2},...,s_n\}$ is graphic and $\sum\limits_{i=1}^{n}s_i \geq 2(n-1)$.  
\end{thrm}

\begin{pf}  ($\Rightarrow$) Suppose that $s$ is connected and graphic. It is required to show that $s'$ is graphic and that $\sum\limits_{i=1}^{n}s_i \geq 2(n-1)$. \\
As $s=d(G)$ for some simple (connected) graph $G$ then $s_i\leq n-1$ for all $i=1,...,n$ and there exists a graph $G'$ with vertex set $V_{G'}=V_G\setminus \{v_1\}$ and edge set $E_{G'}=E_G\setminus \{v_1v_i\mid v_1v_i\in E_G\}$ such that $d(G')=s'$. As $G$ is a simple graph then it follows that $G'$ is also a simple graph, hence $s'$ is graphic. As $s=d(G)$ for some (simple) connected graph $G$, where $|V_G|=n$, then as $G$ is connected $|E_G|\geq n-1$, hence $\sum\limits_{i=1}^{n}s_i \geq 2(n-1)$.

($\Leftarrow$) Suppose that $s'$ is graphic and that $\sum\limits_{i=1}^{n}s_i \geq 2(n-1)$. It is required to show that $s$ is both graphic and connected. \\

As $s'$ is graphic then $s'_i\leq n-2$ for all $i\in\{1,...,n-1\}$. Adding a term $s_1$ (which is necessarily $\leq n-1$) results in $s$ also being graphic as in the worst case scenario i.e. where $s_1=s_n=n-1$, then $s=\{s_1,...,s_n\}=\{n-1,...,n-1\}$ which is the degree sequence of the simple graph $K_{n}$. Suppose that $\sum\limits_{i=1}^{n}s_i = 2(n-1)$, then $s$ is the degree sequence of a graph $G$ with $n-1$ edges and $n$ vertices which means that $G$ is a tree, hence $s$ is connected. If $\sum\limits_{i=1}^{n}s_i > 2(n-1)$ then either $s=d(G)$ for some connected graph $G$ or it is possible to apply Lemma \ref{LH} repeatedly until a graph $G$ is found such that $G$ is connected and $d(G)=s$.
\end{pf}


The following result is what can be thought of as a {\it connected version} of Theorem \ref{HH}. However, note that it is not possible to simply add the word {\it connected} to the statement of Theorem \ref{HH} as $s=(2,2,1,1)$ is connected but $s'=(1,0,1)$ is not connected. 
 
\begin{thrm}\label{Main} Given a sequence $s=\{s_1,...,s_n\}$ of positive integers such that $s_i\geq s_{i+1}$ for $i=1,...,n-1$ then $s$ is connected and graphic if and only if the sequence $s'=\{s'_1,...,s'_{n-1}\}=\{s_1-1,s_2-1,...,s_{s_n}-1, s_{s_n+1},...,s_{n-1}\}$ is connected and graphic.
\end{thrm}


\begin{pf} ($\Rightarrow$) Suppose that $s=\{s_1,...,s_n\}$ is connected. It is required to show that $s'$ is both graphic and connected. \\
To show that $s'$ is graphic it is required to show that $ \sum\limits_{i=1}^{n-1}s'_i $ is even and that all vertices have degree less than or equal to $n-2$. Observe that \[\sum\limits_{i=1}^{n-1}s'_i  = \sum\limits_{i=1}^{n-1}s_i - s_n =  \sum\limits_{i=1}^{n}s_i - 2s_n. \]
As $s$ is graphic then $ \sum\limits_{i=1}^{n}s_i$ is even and so $\sum\limits_{i=1}^{n}s_i - 2s_n$ is also even. As $s$ is graphic then all vertices $v_i\in V_G$ with $i\in\{1,...,n\}$ must satisfy $deg(v_i)\leq n-1$. All vertices with degree $n-1$ in $G$ are necessarily connected to $v_n$ whereas vertices with degree less than $n-1$ may or may not be connected to $v_n$. It follows that after deleting $v_n$ and all edges containing $v_n$ that the maximum degree which any vertex can have in any $G'$ is $n-2$ (where $d(G')=s'$). \\

To show that $s'$ is connected it is required to show that $ \sum\limits_{i=1}^{n-1}s'_i \geq 2(n-2)$ i.e. there exists a graph $G'$ with $d(G')=s'$ and $|V_{G'}|=n-1$ such that $|E_{G'}|\geq n-2$. As $s$ is graphic then $1\leq s_i \leq n-1$.\\

Let $s_{n}=1$: As $s_n=1=\delta_G$ then $s'=(s_1-1, s_2, ..., s_{n-1})$. Not all $s_i=1$ except in the case where $s=(1,1)$ resulting in  $s'=(0)$ which is a connected degree sequence. As $s$ is connected and $deg(v_n)=1$ then $v_n$ is a leaf of a connected graph $G$ and so deleting $v_n$ cannot result in a disconnected graph $G'$, hence $s'$ is connected when $s_n=1$.  \\

Let $s_{n}=k$ where $2\leq k\leq n-1$: As $s_n=k=\delta_G$ then \[s'=(s_1-1, s_2-1,...,s_{s_n}-1, s_{s_n+1} ,..., s_{n-1}).\] Assuming the worst case scenario i.e. $\Delta_G=\delta_G =k$, then this gives \[s'=(\underbrace{k-1,...,k-1}_{k},\underbrace{k,...,k}_{n-k-1})\] which means that \[ \sum\limits_{i=1}^{n-1}s'_i \geq k(k-1) + k(n-k-1)=k(n-2)\geq 2(n-2)\] whenever $2\leq k\leq n-1$. Hence, $s'$ is connected when $s_n=k$ where $2\leq k\leq n-1$.\\


($\Leftarrow$) Suppose that $s'=\{s'_1,...,s'_{n-1}\}=\{s_1-1,s_2-1,...,s_{s_n}-1, s_{s_n+1},...,s_{n-1}\}$ is connected. It is required to show that $s$ is both graphic and connected. \\

 As $s'$ is connected then there exists some $G'$ with $|V_{G'}|=n-1$ and $d(G')=s'$ where the degree of all vertices in $V_{G'}$ is less than or equal to $n-2$. As $V_{G}=V_{G'}\cup \{v_n\}$ (and all $s_i$ are necessarily $\leq n-1$) then all vertices in $V_{G}$ will have degree at most $n-1$. Observe that \[\sum\limits_{i=1}^{n}s_i  =  \sum\limits_{i=1}^{n-1}s'_i + 2s_n. \]
As $s'$ is graphic then $ \sum\limits_{i=1}^{n-1}s'_i $ is even and so $\sum\limits_{i=1}^{n-1}s'_i + 2s_n$ is also even.  \\

As $s'$ is connected then there exists some $G'$ with $d(G')=s'$ where $|E_{G'}|\geq n-2$ as $|V_{G'}|=n-1$. As $s_n\geq 1$ then this means that there is at least one edge in $G$ which has $v_n$ as an endpoint and some $v_{i}$ with $i\in\{1,...,n-1\}$ as the other endpoint.

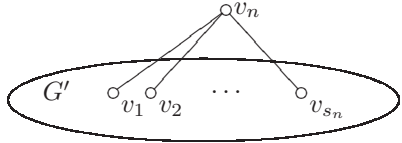
\begin{figure}[H]
\begin{center}
\scalebox{1}{$\begin{xy}\POS (-5,5) *\cir<2pt>{} ="a" *+!UL{v_1},
(0,5) *\cir<2pt>{} ="b" *+!UL{v_2},
(20,5) *\cir<2pt>{} ="c" *+!UL{v_{s_n}},
(10,16) *\cir<2pt>{} ="d" *+!L{v_n},
(-12.5,4)*+!{G'},
 (10,5)*+!{\dots},

\POS "a" \ar@{-}  "d",
\POS "b" \ar@{-}  "d",
\POS "c" \ar@{-}  "d",

\POS(7,4),  {\ellipse(26,5.5)<>{}},

 \end{xy}$}
 
\caption{If $s'$ is connected and graphic then $s$ is connected and graphic }
\label{Proof}
\end{center}
\end{figure}

 This observation along with the fact that $|E_{G}|\geq |E_{G'}|+1>n-2$, where $|V_{G}|=n$, means that $s$ is connected. 
\end{pf}

\begin{exa} An example, using Theorem \ref{Main}, of what Hakimi would term a ``set of successive reduction cycles" is shown in Figure \ref{RedCyc}. 

\begin{figure}[h]
\begin{center}
\scalebox{0.85}{$\begin{xy}
 \POS (10,20) *\cir<2pt>{} ="a" *+!D{} ,
 (20,20) *\cir<2pt>{} ="b" *+!R{},
 (25,10) *\cir<2pt>{} ="c" *+!L{},
 (20,0) *\cir<2pt>{} ="d" *+!L{},
 (10,0) *\cir<2pt>{} ="e" *+!L{},
 (5,10) *\cir<2pt>{} ="f" *+!L{},
 (15,-5) *+!{(4,4,3,3,3,3,3)},
 (33,9) *+!{\longrightarrow}
 
\POS "a" \ar@{-}  "b",
\POS "a" \ar@{-}  "c",
\POS "a" \ar@{-}  "d",
\POS "a" \ar@{-}  "e",
\POS "b" \ar@{-}  "d",
\POS "b" \ar@{-}  "e",
\POS "b" \ar@{-}  "f",
\POS "c" \ar@{-}  "e",
\POS "c" \ar@{-}  "f",
\POS "d" \ar@{-}  "f",

\POS (40,20) *\cir<2pt>{} ="a1" *+!D{} ,
 (50,20) *\cir<2pt>{} ="b1" *+!R{},
 (50,10) *\cir<2pt>{} ="c1" *+!L{},
 (45,2.5) *\cir<2pt>{} ="d1" *+!L{},
 (40,10) *\cir<2pt>{} ="e1" *+!L{},
 (45,-5) *+!{(3,3,3,3,2)},
 (58,9) *+!{\longrightarrow}
 
\POS "a1" \ar@{-}  "b1",
\POS "a1" \ar@{-}  "c1",
\POS "a1" \ar@{-}  "e1",
\POS "b1" \ar@{-}  "c1",
\POS "b1" \ar@{-}  "e1",
\POS "c1" \ar@{-}  "d1",
\POS "d1" \ar@{-}  "e1",

\POS  (65,15) *\cir<2pt>{} ="a2" *+!D{} ,
 (75,15) *\cir<2pt>{} ="b2" *+!R{},
 (75,5) *\cir<2pt>{} ="c2" *+!L{},
 (65,5) *\cir<2pt>{} ="d2" *+!L{},
 (70,-5) *+!{(3,3,2,2)},
 (83,9) *+!{\longrightarrow}
 
\POS "a2" \ar@{-}  "b2",
\POS "a2" \ar@{-}  "d2",
\POS "b2" \ar@{-}  "c2",
\POS "b2" \ar@{-}  "d2",
\POS "c2" \ar@{-}  "d2",

\POS (90,15) *\cir<2pt>{} ="a3" *+!D{} ,
 (100,15) *\cir<2pt>{} ="b3" *+!R{},
 (90,5) *\cir<2pt>{} ="c3" *+!L{},
 (95,-5) *+!{(2,2,2)},
 (107,9) *+!{\longrightarrow},

\POS "a3" \ar@{-}  "b3",
\POS "a3" \ar@{-}  "c3",
\POS "b3" \ar@{-}  "c3",

\POS (115,10) *\cir<2pt>{} ="a4" *+!D{} ,
 (125,10) *\cir<2pt>{} ="b4" *+!R{},
 (120,-5) *+!{(1,1)},
  (132,9) *+!{\longrightarrow},
\POS "a4" \ar@{-}  "b4",

\POS (140,10) *\cir<2pt>{} ="a4" *+!D{} ,
(140,-5) *+!{(0)},

 \end{xy}$ }
 \caption{A set of successive reduction cycles}
\label{RedCyc}
\end{center}
\end{figure}
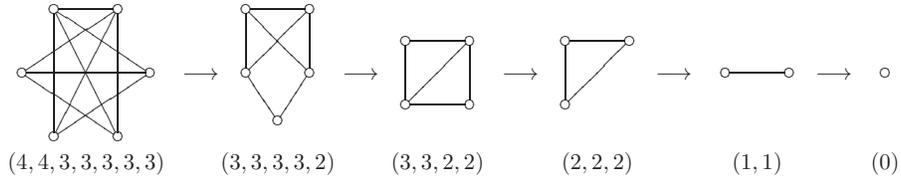
 
\end{exa} 

\section{Comments}

Theorem \ref{HH} is used in \cite{Hk} to check algorithmically when a given sequence $s$ is graphic. In a similar manner, Theorem \ref{Main} suggests an algorithm which can be used to determine if a given sequence $s$ is connected and graphic.


\bibliographystyle{plain}      
\bibliography{refs2015}

\end{document}